\documentclass[12pt]{article}

\usepackage{amssymb}
\usepackage{latexsym}

\title{ On Zumkeller Numbers}
\author{
K.P.S. Bhaskara Rao \ \ \ \ \ \ \ \ Yuejian Peng \\
Department of Mathematics and Computer Science \\
Indiana State University \\
Terre Haute, IN, 47809, USA \\
Email: bkopparty@indstate.edu;  yuejian.peng@indstate.edu }
\date{}

\newtheorem{defi}{Definition}
\newtheorem{theo}{Theorem}

\newtheorem{remark}[theo]{Remark}
\newtheorem{lemma}[theo]{Lemma}

\newtheorem{coro}[theo]{Corollary}

\newtheorem{prop}[theo]{Proposition}
\newtheorem{fact}[theo]{Fact}

\newcommand{\qed}{\hspace*{\fill} \rule{7pt}{7pt}}

\begin{document}
\maketitle
\begin{abstract}
Generalizing the concept of a perfect number, Sloane's sequences of
integers A083207 lists the sequence of integers $n$ with the property: the
positive factors of $n$ can be partitioned into two disjoint parts so that the
sums of the two parts  are equal. Following [4] Clark et al., we shall call such integers, Zumkeller numbers. Generalizing
this, in [4] a number $n$ is called a half-Zumkeller number if the positive proper factors of $n$ can be partitioned into two disjoint parts so that the
sums of the two parts  are equal.

An extensive study of properties of Zumkeller numbers, half-Zumkeller numbers and their relation
to practical numbers is undertaken in this paper.

In [4] Clark et al., announced results about Zumkellers numbers  and half-Zumkeller numbers 
and suggested two conjectures. In the present paper we shall settle one of the conjectures, prove the second conjecture in some special cases and
prove several results related to the second conjecture. We shall also show that if there
is an even Zumkeller number that is not half-Zumkeller it should be bigger than $7.2334989 \times 10^9$.

\end{abstract}

\section{Introduction}

A positive integer $n$ is called a {\em perfect} number if $n$ equals the sum of its proper
positive factors. 

Generalizing this concept in 2003, Zumkeller   published in Sloane's
sequences of integers A083207
a sequence of integers n with the property that the positive
factors of n can be partitioned into two disjoint parts so that the
sums of the two parts are equal.

Following Clark et al., we shall call such integers, Zumkeller numbers.
In \cite{walsh} Clark et al., announced several results about Zumkeller
numbers and suggested some conjectures in order to understand Zumkeller numbers.

In section 2 of this paper we shall study some properties of Zumkeller numbers.
In section 3 we shall study the relations between practical numbers and Zumkeller
numbers. We shall also settle a conjecture from [4] in this section. In
section 4 we shall define and study the properties of
half-Zumkeller numbers and prove results analogous to results about
Zumkeller numbers for half-Zumkeller numbers. We shall also make substantial
contributions
regarding the second conjecture from [4].

\section{Zumkeller numbers}
\begin{defi}
A positive integer $n$ is said to be a {\em Zumkeller} number if the positive
factors of n can be partitioned into two disjoint parts so that the
sums of the two parts are equal. We shall call such a partition a {\em Zumkeller
partition}.
\end{defi}

6, 12, 20, 24, 28, 30, 40 are the first few Zumkeller numbers.
We shall start with a few simple facts. We shall also designate some of the
known results as facts and for some of the facts we shall provide the proofs.

Let $\sigma(n)$ represent the sum of all positive factors of $n$.

\begin{fact}\label{sigma}
Let the prime factorization of  $n$ be $\Pi_{i=1}^m p_i^{k_i}$. Then 
$$\sigma(n)=\Pi_{i=1}^m (\sum_{j=0}^{k_i} p_i^j)$$ and
$${\sigma(n)\over n}=\Pi_{i=1}^m (\sum_{j=0}^{k_i} p_i^{j-k_i})=\Pi_{i=1}^m (\sum_{j=0}^{k_i} p_i^{-j})\le \Pi_{i=1}^m {p_i \over p_i-1}.$$ 
\end{fact}

$\sigma(n)$ of a Zumkeller number $n$ cannot be odd. This is observed in the next fact.

\begin{fact}\label{f3}\cite{walsh}
If $n$ is a Zumkeller number, then $\sigma(n)$ must be even. Therefore the number of positive odd factors of $n$ must be even.
\end{fact}
 
{\bf Proof.} If $\sigma(n)$ is odd, then it is impossible to partition the positive factors of $n$ into two equal-summed parts. If the number of positive odd factors of $n$ is odd, then $\sigma(n)$ is odd. \qed
 
The following fact gives a necessary and sufficient condition for $n$ to be a Zumkeller number.

\begin{fact}\label{f1}
 $n$ is a Zumkeller  number if and only if ${\sigma(n)-2n \over 2}$ is
either zero or is a sum of distinct positive  factors of $n$ excluding  $n$ itself.
\end{fact}

{\bf Proof.} $n$ is a Zumkeller  number if and only if there exists $a$ which is either $0$ or a sum of some factors of $n$ excluding  $n$ itself such that
\begin{equation}\label{eq1}
n+a=\sigma(n)-(n+a).
\end{equation}
This is equivalent to ${\sigma(n)-2n \over 2}=a$. \qed

Let us state  a simple necessary condition for a number to be Zumkeller. This says that every Zumkeller number is abundant.

\begin{fact}\label{f2}\cite{walsh}
If $n$ is a Zumkeller number, then $\sigma(n)\ge 2n$.
\end{fact}

This follows from Fact \ref{f1}.

Based on Fact \ref{f3}, we shall now give a necessary condition
for an even number to be Zumkeller.

\begin{fact}\label{f4}
Let  the prime factorization of an even Zumkeller number $n$ be $2^kp_1^{k_1}p_2^{k_2}\cdots p_m^{k_m}$. Then at least one of $k_i$ must be odd.
\end{fact}

{\bf Proof.} Note that the number of positive odd factors of $n$ is $(k_1+1)(k_2+1)\cdots (k_m+1)$. At least one of $k_i$ must be odd in order to make the number $(k_1+1)(k_2+1)\cdots (k_m+1)$ be even. \qed

The following fact from \cite{walsh} gives a method of generating new Zumkeller numbers from a known Zumkeller number.

\begin{fact}\cite{walsh}\label{fact55} If $n$ is a Zumkeller number and $p$ is a prime with $(n, p)=1$, then $np^l$ is Zumkeller for any positive integer $l$. 

\end{fact}

{\bf Proof.} Since $n$ is a Zumkeller number, the set of all positive factors of $n$, denoted by $D_0$  can be partitioned  into two disjoint parts $A_0$ and $B_0$ so that the sums of the two parts are equal. 
Group the factors of $np^l$ into $l+1$ groups $D_0$, $D_1$, $\ldots$, $D_{l}$  according to how many factors of $p$ they admit, i.e., $D_i$ consists of all positive  factors of $np^l$ admitting $i$ factors of $p$ for every $i$, $0\le i\le l$. Then each $D_i$ can be partitioned into two disjoint parts so that the sum of these two parts  equal according to the Zumkeller partition of $D_0=A_0\cup B_0$. \qed

We give another method of  generating new Zumkeller numbers from a known Zumkeller number.

\begin{prop} \label{prop6} Let $n$ be a Zumkeller number and
$p_1^{k_1}p_2^{k_2}\cdots p_m^{k_m}$ be the prime factorization of  $n$.
Then  for any positive integers $l_1, \ldots, l_m$,
$$p_1^{k_1+l_1(k_1+1)}p_2^{k_2+l_2(k_2+1)}\cdots p_m^{k_m+l_m(k_m+1)}$$ is
Zumkeller.
\end{prop}

{\bf Proof.} It is sufficient to show that $p_1^{k_1+l_1(k_1+1)}p_2^{k_2}\cdots p_m^{k_m}$ is Zumkeller
if $n=p_1^{k_1}p_2^{k_2}\cdots p_m^{k_m}$ is Zumkeller.
Assume that $n=p_1^{k_1}p_2^{k_2}\cdots p_m^{k_m}$ is Zumkeller.
Then the set of all positive factors of $n$, denoted by $D_0$ can be
partitioned into two parts $A_0$ and $B_0$ so that the sums of
these two parts  are equal. Note that the factors of
$p_1^{k_1+l_1(k_1+1)}p_2^{k_2}\cdots p_m^{k_m}$ can be partitioned
into $l_1+1$ disjoint groups $D_i$, $0\le i\le l_1$, where
elements in $D_i$, $1\le i\le l_1$ are obtained by multiplying elements in $D_0$ with
$p_1^{i(k_1+1)}$. Note that every element in
$D_i$ admits at least $i(k_1+1)$ and at most $k_1+i(k_1+1)$ factors
of $p_1$. So for $i \ne j$,  $D_i$ and $D_j$ are disjoint.
Since $D_0$ can be partitioned into two disjoint parts $A_0$ and $B_0$
so that the sums of  the two parts  are equal, all the  other $D_i$'s
can be partitioned into two disjoint parts $A_i$ and $B_i$ correspondingly
so that the sums of  these two parts   equal .
 This proves that $p_1^{k_1+l_1(k_1+1)}p_2^{k_2}\cdots p_m^{k_m}$ is a
 Zumkeller number. \qed

\section{Practical numbers and Zumkeller numbers}
Fact \ref{f1} states that if a certain number related to $n$ is a sum of
distinct positive factors of $n$, then $n$ is Zumkeller.
From the definition of practical numbers (given below),
it is natural to consider the relations between practical numbers and
Zumkeller numbers. We shall do this in this section.
Practical numbers were introduced in \cite{Srinivasan}.

\begin{defi}
A  positive integer $n$ is said to be a {\em practical} number if all  positive integers less than $n$ can be represented as sums of distinct factors of $n$. 
\end{defi}

As Stewart (1954) in \cite{stewart} showed , it is straightforward to determine whether a number is practical from its prime factorization.
\begin{fact}\label{pract1}\cite{stewart}
A positive integer $n$ with the prime factorization
$p_1^{k_1}p_2^{k_2}\cdots p_m^{k_m}$ and $p_1<p_2<\ldots<p_m$ is
a practical number  if and only if $p_1 = 2$ and
for $1\le i\le m-1$, $p_{i+1}\le \sigma(p_1^{k_1}\cdots p_i^{k_i})+1$.
\end{fact}

The following result from \cite{rao} and \cite{stewart} says that every positive integer $ \le n$ is a sum of distinct
factors on $n$ if and only if every positive integer $ \le \sigma(n)$ is a sum of
distinct factors of $n$.


\begin{prop}\cite{rao}\label{very}
A positive integer $n$ is a practical number if and only if every positive integer $ \le \sigma(n)$ is a sum of distinct factors of $n$.
\end{prop}

 Note that every practical number is even.
 Also note that for every positive integer $k$, $2^{k}$ is practical.
 The next result gives a sufficient condition for the product of $2^k$ and a prime
 $p \ne 2$ to be a Zumkeller number.

\begin{fact}\cite{walsh}\label{walshpractical}
For any prime $p\neq 2$ and positive integer $k$ with $p\le 2^{k+1}-1$, $2^kp$ is a Zumkeller number.
\end{fact}

With the aim of generalizing the above result for all practical
numbers, the following problem was proposed as a conjecture in \cite{walsh}.

\bigskip

{\bf Conjecture 1} \cite{walsh} \ 
Let $n$ be a practical number and $p$ be a prime with $(n, p)=1$,  and
$p\le \sigma(n)$. Then, $np$ is Zumkeller.

Clark et al., suggested that some further restrictions on $\sigma(n)$ or
${\sigma(n) \over n}$ are possibly needed for the conjecture to be true.

\bigskip

We shall prove a comprehensive theorem (Theorem \ref{theopz}) that  settles this conjecture.
We shall first
find necessary and sufficient conditions for a practical number to be a Zumkeller
number.

\begin{prop}\label{proppraczu}
A practical number $n$ is a Zumkeller number if and only if $\sigma(n)$ is even.

\end{prop}
{\bf Proof.} If $n$ is Zumkeller, then $\sigma(n)$ is even by Fact \ref{f3}. 

If $\sigma(n)$ is even, then ${\sigma(n) \over 2}$ is a positive integer smaller than $\sigma(n)$. Since $n$ is practical, by Proposition \ref{very},  ${\sigma(n) \over 2}$ is
the sum of some positive factors of $n$. Therefore, the set of all
positive factors of $n$ can be partitioned into two equal-summed parts
and $n$ is Zumkeller. \qed

\begin{theo}\label{theopz}
Let $n$ be a practical number and $p$ be a prime with $(n, p)=1$.
Let $l$ be a positive integer.  Then,

(i) If  $\sigma(n)$ is even, then $np^l$ is Zumkeller.

(ii) If  $\sigma(n)$ is odd,   then $np^l$ is Zumkeller if and only if
$p\le\sigma(n)$ and $l$ is odd.
\end{theo}

{\bf Proof.}  (i) By Proposition \ref{proppraczu} and Fact \ref{fact55}, $np^l$is Zumkeller.

(ii) Let $\sigma(n)$ be odd. We first apply the conclusion in Case 1 to show that $np^l$ is Zumkeller
if $p\le\sigma(n)$ and $l$ is odd. Since $n$ is practical and
$p \le \sigma(n)$, by Fact \ref{pract1}, $np^l$ is practical.
Since $l$ is odd and $p$ is odd, $\sum_{i=0}^l p^i$ is even and consequently,
$\sigma(np^l)=(\sum_{i=0}^l p^i) \sigma(n)$ is even.  By Proposition
\ref{proppraczu}, we obtain that $np^l$ is Zumkeller.

Now we shall show that if $np^l$ is Zumkeller then,  $p\le\sigma(n)$ and $l$ is odd.
Let $np^l$ be Zumkeller. Then for every integer $i$,  $0\le i\le l$,  there exist $a_i, b_i$ each of which is either $0$ or a sum of some factors of $n$ such that
$$\sum_{i=0}^l p^i(a_i-b_i)=0;$$
and 
$$a_i+b_i=\sigma(n).$$
Therefore, $p$ divides $a_0-b_0$. Since $a_0+b_0=\sigma(n)$ is odd, $a_0-b_0\neq 0$. Therefore, $\vert a_0-b_0\vert \ge p$. Hence, $$\sigma(n)=a_0+b_0\ge \vert a_0-b_0\vert \ge p.$$ 
Since $np^l$ is Zumkeller, $\sigma(np^l)=(\sum_{i=0}^l p^i) \sigma(n)$ must be even by Fact \ref{f3}. Since $\sigma(n)$ is odd, $\sum_{i=0}^l p^i$ must be even.
Therefore, $l$ must be odd. \qed

In Proposition \ref{proppraczu} we showed that every practical number $n$ with even $\sigma(n)$ is Zumkeller. But,

\begin{remark} 
There are Zumkeller numbers  that are not practical numbers.
\end{remark}

All odd Zumkeller numbers  are not practical numbers.
For example,
$${\sigma(945)-2\cdot 945 \over 2}={(1+3+3^2+3^3)\cdot (1+5)\cdot (1+7)-1890 \over 2}=15$$ is a
factor of 945.   By Fact \ref{f1}, 945 is a Zumkeller number.
But 945 is not a practical number since all practical numbers are even.

Also, $70=2\cdot5\cdot7$ is an even Zumkeller number that is not practical.
In fact,
$${\sigma(70)-2\cdot 70 \over 2}={(1+2)\cdot (1+5)\cdot (1+7)-140 \over 2}={144-140 \over 2}= 2$$ 
 is a  factor of 70. By Fact \ref{f1},
 70 is a Zumkeller number. But it is not a practical number by
 Fact \ref{pract1}.
 
\bigskip

Now we shall develop several results for Zumkeller numbers and these will be used later.
A finer analysis of the proof of (ii) in the above theorem,
gives us the following result.

\begin{prop}\label{prop4}
Let $l$ be a positive integer. Let $n$ be  a non-Zumkeller
number and $p$ be a prime with $(n, p)=1$.  If $np^l$ is Zumkeller,
then $p\le \sigma(n)$.  If $np^l$ is Zumkeller and $\sigma(n)$ is odd,
then $l$ is odd.
\end{prop}

{\bf Proof.} This is similar to the proof of the necessity part of
Theorem \ref{theopz} (ii). The only modification is to replace
`Since $a_0+b_0=\sigma(n)$ is odd, $a_0-b_0\neq 0$' by `Since $a_0+b_0=\sigma(n)$ and $n$ is not Zumkeller, $a_0-b_0\neq 0$'. \qed

We will now give several necessary and sufficient conditions for $np$ to be a Zumkeller number.

\begin{prop}\label{multiplyz}
Let $n$ be  a positive integer and $p$ be a prime with $(n, p)=1$.  Then the following conditions are equivalent:

(i) $np$ is Zumkeller.

(ii) The set  of all positive factors of $n$ can be partitioned into two
disjoint parts  $D_1$ and $D_2$ such that $p(\sum_{d\in D_2} d - \sum_{d\in D_1} d)$ can be
written as a sum of some factors of $n$   minus  the sum of the rest of the
factors of $n$.

(iii) The set  of all positive factors of $n$ can be partitioned into two
disjoint parts  $D_1$ and $D_2$ such that ${(p+1)(\sum_{d\in D_2} d - \sum_{d\in D_1} d) \over 2}$ can be
written as a sum of some elements in $D_2$ minus  a sum of some
elements in $D_1$.

(iv) The set  of all positive factors of $n$ can be partitioned into 4 disjoint parts  $A_1$, $A_2$,  $A_3$ and $A_4$ such that $(p+1)\sum_{d\in A_1} d +(p-1)\sum_{d\in A_2} d=(p+1)\sum_{d\in A_3} d +(p-1)\sum_{d\in A_4} d$.
\end{prop}

{\bf Proof.}  $np$ is Zumkeller if and only if the set $D_0$ of all positive
factors of $n$ can be partitioned into $D_0=D_1 \cup D_2$ and $D_0=D_3 \cup D_4$ such that

$$p( \sum_{d\in D_1} d) + \sum_{d\in D_4} d   = p(\sum_{d\in D_2} d) + \sum_{d\in D_3} d.$$

This can be rewritten as
$$p(\sum_{d\in D_2} d - \sum_{d\in D_1} d)+(\sum_{d\in D_3} d - \sum_{d\in D_4} d)=0.$$
This is equivalent to 
\begin{eqnarray}\label{eq4}
p(\sum_{d\in D_2} d - \sum_{d\in D_1} d)=-(\sum_{d\in D_3} d - \sum_{d\in D_4} d).
\end{eqnarray}
This is equivalent to (ii).

Now, by adding $\sum_{d\in D_2} d - \sum_{d\in D_1} d$ to both sides of the above equation
we get,

\begin{eqnarray*}
p(\sum_{d\in D_2} d - \sum_{d\in D_1} d)+(\sum_{d\in D_2} d - \sum_{d\in D_1} d)&=&(\sum_{d\in D_2} d - \sum_{d\in D_1} d)-(\sum_{d\in D_3} d - \sum_{d\in D_4} d)\\
&=&(\sum_{d\in D_2} d -\sum_{d\in D_3} d)-(\sum_{d\in D_1} d-\sum_{d\in D_4} d)\\
&=&(\sum_{d\in D_2} d-\sum_{d\in D_3\cap D_2} d)-\sum_{d\in D_3\cap D_1} d\\
&&-(\sum_{d\in D_1} d-\sum_{d\in D_4\cap D_1} d)+\sum_{d\in D_4\cap D_2} d \\
&=&\sum_{d\in D_4\cap D_2} d-\sum_{d\in D_3\cap D_1} d-\sum_{d\in D_3\cap D_1} d+\sum_{d\in D_4\cap D_2} d \\
&=&2(\sum_{d\in D_4\cap D_2} d-\sum_{d\in D_3\cap D_1} d).
\end{eqnarray*}

Thus,
$${(p+1)(\sum_{d\in D_2} d - \sum_{d\in D_1} d) \over 2}=\sum_{d\in D_4\cap D_2} d-\sum_{d\in D_3\cap D_1} d.$$
This is equivalent to (iii). To go in the converse direction one needs to take
$D_3 = (D_3 \cap D_1)\cup (D_2 - D_4 \cap D_2)$ and
$D_4 = (D_4 \cap D_2)\cup (D_1 - D_3 \cap D_1)$.

By adding $p(\sum_{d\in D_3} d - \sum_{d\in D_4} d)$    to both sides of equation (\ref{eq4})
we get that equation (\ref{eq4}) is equivalent to 

$$p(\sum_{d\in D_2} d - \sum_{d\in D_1} d)+p(\sum_{d\in D_3} d - \sum_{d\in D_4} d)=-(\sum_{d\in D_3} d - \sum_{d\in D_4} d)+p(\sum_{d\in D_3} d - \sum_{d\in D_4} d).$$
This is equivalent to 
$$p(\sum_{d\in D_2} d - \sum_{d\in D_4} d+\sum_{d\in D_3} d - \sum_{d\in D_1} d)=(p-1)(\sum_{d\in D_3} d - \sum_{d\in D_4} d).$$
This is equivalent to 
$$p(\sum_{d\in D_2} d - \sum_{d\in D_2\cap D_4} d -\sum_{d\in D_1\cap D_4} d+\sum_{d\in D_3} d - \sum_{d\in D_1\cap D_3} d- \sum_{d\in D_1\cap D_4} d)=(p-1)(\sum_{d\in D_3} d - \sum_{d\in D_4} d).$$
This is equivalent to 
$$p(\sum_{d\in D_2\cap D_3} d -\sum_{d\in D_1\cap D_4} d+ \sum_{d\in D_2\cap D_3} d- \sum_{d\in D_1\cap D_4} d)=(p-1)(\sum_{d\in D_3} d - \sum_{d\in D_4} d).$$
This is equivalent to 
$$2p(\sum_{d\in D_2\cap D_3} d - \sum_{d\in D_1\cap D_4} d)=(p-1)(\sum_{d\in D_3} d - \sum_{d\in D_4} d).$$
This is equivalent to 
\begin{eqnarray*}
(p+1)(\sum_{d\in D_2\cap D_3} d - \sum_{d\in D_1\cap D_4} d)&=&(p-1)(\sum_{d\in D_3} d - \sum_{d\in D_4} d-\sum_{d\in D_2\cap D_3} d +\sum_{d\in D_1\cap D_4} d) \\
&=&(p-1)(\sum_{d\in D_1\cap D_3} d - \sum_{d\in D_2\cap D_4} d).
\end{eqnarray*}
This is equivalent to 
$$(p+1)\sum_{d\in D_2\cap D_3} d+(p-1)\sum_{d\in D_2\cap D_4} d=(p+1)\sum_{d\in D_1\cap D_4} d+(p-1)\sum_{d\in D_1\cap D_3} d.$$
By letting $A_1=D_2\cap D_3$, $A_2=D_2\cap D_4$, $A_3=D_1\cap D_4$, and $A_4=D_1\cap D_3$ and observing that
$\{A_1, A_2, A_3, A_4 \}$ is a partition of the positive factors of $n$
we see that the above equation is equivalent to (iv).

\qed

\bigskip

\begin{remark} If the set  of all positive factors of $n$ can be partitioned into two disjoint parts  $D_1$ and $D_2$ such that  $\sum_{d\in D_2} d - \sum_{d\in D_1} d$ is small, then it is usually easy to check  condition (iii) in Proposition \ref{multiplyz}.
\end{remark}
Let us see an application of Proposition \ref{multiplyz}.

\begin{fact}\label{257}
 $2\times 5^2\times 7^2\times p$ is Zumkeller for $p=11, 13$. 
\end{fact}

{\bf Proof.} The set of all positive factors of $2\times 5^2\times 7^2$ can be partitioned into two disjoint parts $D_1=\{2450, 98, 50, 35, 10, 5, 2\}$ and $D_2=\{1225, 490, 350, 245, 175, 70, 49, 25, 14, 7, 1\}$ and $\sum_{d\in D_2} d - \sum_{d\in D_1} d=1$.

${(11+1)(\sum_{d\in D_2} d - \sum_{d\in D_1} d) \over 2}=6=(7+1)-2$ is a sum of some elements in $D_2$  minus  a sum of some elements in $D_1$ and it satisfies Proposition \ref{multiplyz}(iii) for $p=11$. So $2\times 5^2\times 7^2\times 11$ is Zumkeller.

${(13+1)(\sum_{d\in D_2} d - \sum_{d\in D_1} d) \over 2}=7$ is a sum of some elements in $D_2$  and it satisfies Proposition \ref{multiplyz}(iii) for $p=13$. So $2\times 5^2\times 7^2\times 13$ is Zumkeller. \qed




One can use the above method to show that Fact \ref{257} holds for more prime numbers.

\bigskip

In Proposition 1.4 of \cite{rao} a condition on the existence of factors
with certain nice properties was shown to be sufficient for an integer to be a practical number.
We
shall prove a similar theorem for Zumkeller numbers. In this case we shall also develop a method of finding a
Zumkeller partition.

\begin{prop}\label{prop7}
If $1=a_1<a_2<\cdots<a_k=n$ are all factors of $n$ with $a_{i+1}\le 2a_i$ for all $i$ and $\sigma(n)$ is even, then $n$ is Zumkeller. 
\end{prop}

{\bf Proof.} Starting with a positive sign for $a_k=n$, we will assign positive or negative signs to each $a_i$ and  show that the sum of all $a_i$, $1\le i\le k$  with the positive  or negative signs assigned is 0.
Then it will imply that $\sigma(n)$ can be partitioned into two equal-summed parts.

$a_k$ is assigned positive sign, and $a_{k-1}$ is assigned the negative sign. Note that 
$0<a_k-a_{k-1}\le a_{k-1}$ since $a_{k}\le 2a_{k-1}$. Since the current sum  $a_k-a_{k-1}$  is positive, we assign the negative sign to $a_{k-2}$. Then $-a_{k-2}<a_k-a_{k-1}-a_{k-2}\le a_{k-1}-a_{k-2}\le a_{k-2}$ since 
$a_{k-1}\le 2a_{k-2}$. If $a_k-a_{k-1}-a_{k-2}\ge 0$, we assign the negative sign to $a_{k-3}$. Otherwise we assign the positive sign to $a_{k-3}$.
Let $s_i$ be the current sum up to $a_{i}$. In general, the sign assigned to $a_{i-1}$ is opposite  of the sign of $s_i$.
Let us show inductively that $\vert s_i\vert \le a_{i}$ for $1\le i\le k$. It is true for $i=k$. Assume that  $\vert s_{i+1}\vert \le a_{i+1}$. Since the sign of  $a_{i}$ is opposite  of the sign of
$s_{i+1}$, $\vert s_{i}\vert=\vert\vert s_{i+1}\vert -a_{i}\vert$. Note that $-a_{i}\le \vert s_{i+1}\vert -a_{i}\le a_{i+1}-a_{i}\le a_i$ since $a_{i+1}\le 2a_i$. Therefore $\vert s_{i}\vert\le a_i$.  So $\vert s_{1}\vert\le a_1=1$. Since $\sigma(n)$ is even, $s_1$, which is obtained by assigning a positive or negative sign to each of the terms in $\sigma(n)$ is even as well.
So $s_1=0$. This implies that $\sigma(n)$ can be partitioned into two equal-summed parts, i.e., $n$ is Zumkeller.
\qed

Clearly, if $n$ is an integer for which we can find factors $1=b_1<b_2<\cdots<b_k=n$ with the property  that $b_{i+1} \le 2b_i$ for all $i$ and $\sigma(n)$ is even, then
the hypothesis of Proposition \ref{prop7} is satisfied (for all factors of $n$). Hence such an $n$ is Zumkeller.

\begin{remark}\label{n!z} In \cite{rao}, $n!$ was shown to be a practical number for all
$n\ge 3$. We shall now apply Proposition \ref{prop7} and show that
$n!$ is a Zumkeller number for all $n\ge 3$.
\end{remark}

{\bf Proof.} If $1=a_1<a_2<\cdots<a_k=n!$ are all factors of $n!$, then clearly $a_{i+1}\le 2a_i$ for all $i$ (Note that
$1, 2, 3, \ldots, n, 2n, 3n, \ldots, (n-1)n, 2(n-1)n, 3(n-1)n, \ldots, (n-2)(n-1)n, 
2(n-2)(n-1)n, 3(n-2)(n-1)n, \ldots, \prod_{i=3}^n i, 2\prod_{i=3}^n i=n!$
are some of the factors of $n!$.) By Proposition \ref{prop7}, it is sufficient to show
that $\sigma(n!)$ is even for $n\ge 3$.    Let $p$ be the largest
prime $ \le n$,
If $n!=p_1^{k_1}p_2^{k_2}\cdots p_{m_n}^{k_{m_n}}$ is the prime
factorization of $n!$ with $2= p_1<p_2<\cdots <p_{m}$, then clearly $p_{m} = p$.  If $k_m\ge 2$, then $2p\le n$ and by Bertrand's postulate, there exists a prime number $q$ such that ${n \over 2}<q\le n$ and this contradicts the choice of $p$. So $k_{m}= 1$.  This implies that $\sigma(n!)$ is even.

\qed

\bigskip

We shall now discuss odd Zumkeller numbers.
As stated in \cite{walsh}, odd Zumkeller numbers do exist.
In fact, the first several odd abundant numbers with even $\sigma$-value are
Zumkeller. Using Proposition \ref{prop6}, starting with an odd Zumkeller
number one can produce an infinite sequence of Zumkeller numbers.
The next result describes the prime factors of an odd Zumkeller number with a small number of
prime factors.

\begin{fact}\label{f19}
Let the prime factorization of an odd number $n$ be $p_1^{k_1}p_2^{k_2}\cdots p_m^{k_m}$, where $3\le p_1<p_2<\cdots <p_m$. If $n$  is Zumkeller, then 
$$\Pi_{i=1}^m {p_i \over p_i-1}\ge 2.$$
and $m$ is at least $3$. In particular:

1. If $m\le 6$, then $p_1=3$, $p_2=5$,  $7$ or $11$.

2. If $m\le 4$, then $p_1=3$, $p_2=5$. 

3. If $m=3$, then $p_1=3$, $p_2=5$,  and $p_3=7$ or $11$ or $13$.

\end{fact}

{\bf Proof.} If $n$ is Zumkeller, then by Facts \ref{f2} and \ref{sigma}, 
$$2p_1^{k_1}p_2^{k_2}\cdots p_m^{k_m}=2n\le \sigma (n)=\Pi_{i=1}^m (\sum_{j=0}^{k_i} p_i^j).$$
Dividing both sides by $p_1^{k_1}p_2^{k_2}\cdots p_m^{k_m}$, we get 
$$2\le \Pi_{i=1}^m (\sum_{j=0}^{k_i} p_i^{j-k_i})\le \Pi_{i=1}^m {p_i \over p_i-1}.$$
If $m\le 2$, then 
$$\Pi_{i=1}^m {p_i \over p_i-1}\le {3 \over 2}\times {5 \over 4}<2.$$
So $m\ge 3$.
The truth of $\it{1-3}$ follows by verifying the condition $\Pi_{i=1}^m {p_i \over p_i-1}\ge 2$ directly as given below. 

$1$. Let  $m\le 6$. If $p_1\neq 3$, then $p_1\ge 5$ and   
$$\Pi_{i=1}^m {p_i \over p_i-1}\le {5 \over 4}\times {7 \over 6}\times {11 \over 10}\times {13 \over 12}\times {17 \over 16}\times {19 \over 18}<2.$$
Therefore, $p_1=3$. If $p_2>11$, then   $p_2\ge 13$ and
$$\Pi_{i=1}^m {p_i \over p_i-1}\le {3 \over 2}\times {13 \over 12}\times {17 \over 16}\times {19 \over 18}\times {23 \over 22}\times {29 \over 28}<2.$$
Hence, $p_2\le 11$. This implies that $p_2=5$, $7$ or $11$. 

$2$. Let $m\le 4$. By $1$,  $p_1=3$. 
If $p_2>7$, then   $p_2\ge 11$, so
$$\Pi_{i=1}^m {p_i \over p_i-1}\le {3 \over 2}\times {11 \over 10}\times {13 \over 12}\times {17 \over 16}<2.$$
Therefore, $p_2\le7$. This implies that $p_2=5$ or $7$.

$3$. Let $m=3$. By $1$, $p_1=3$. 
If $p_2\neq 5$, then  $p_2\ge 7$ and $p_3\ge 11$. So
$$\Pi_{i=1}^3 {p_i \over p_i-1}\le {3 \over 2}\times
{7 \over 6}\times {11 \over 10}<2.$$
Hence, $p_2=5$. 

If $p_3\ge 17$, then  
$$\Pi_{i=1}^3 {p_i \over p_i-1}\le {3 \over 2}\times {5 \over 4}\times {17 \over 16}<2.$$
Hence, $p_3<17$ and consequently $p_3=7, 11$ or $13$.   \qed

We do not know too much more about odd Zumkeller numbers.
\section{Half-Zumkeller numbers}
\begin{defi} 
A positive integer $n$ is is said to be {\em half-Zumkeller} number if the
proper
positive factors of $n$ can be partitioned into two disjoint parts so that the sums 
of the two parts are equal.
\end{defi}

We shall start with some simple observations.
\begin{fact}\label{h1}
A positive integer $n$ is  half-Zumkeller if and only if
${\sigma(n)-n \over 2}$ is  the sum of some distinct positive  proper  factors
of $n$.
\end{fact}

\begin{fact}\label{h3}
A positive even integer $n$ is  half-Zumkeller if and only if
${\sigma(n)-2n \over 2}$ is zero or the sum of some distinct factors of
$n$ excluding $n$ and ${n \over 2}$.
\end{fact}

{\bf Proof.} Let $n$ be even.  $n$ is half-Zumkeller if and only if there exists $a$ which is zero
or the sum of some  factors of $n$ excluding $n$ and ${n \over 2}$ such that
$${n \over 2}+a=\sigma(n)-n-({n \over 2}+a).$$
Therefore, $a={\sigma(n)-2n \over 2}$. \qed

\begin{fact}\label{257h1}$2\cdot 5\cdot 7$, $2\cdot 5^2\cdot 7$, and $2\cdot 5\cdot 7^2$ are half-Zumkeller numbers.
\end{fact}
 
{\bf Proof.}
$2\cdot 5\cdot 7$ is half-Zumkeller since 
$${\sigma(2\cdot 5\cdot 7)-2(2\cdot 5\cdot 7) \over 2}=2$$ is a sum of factors of
$2\cdot 5\cdot 7$ excluding $2\cdot 5\cdot 7$ and $5\cdot 7$ (by Fact \ref{h3}). Similarly, $2\cdot 5^2\cdot 7$ is half-Zumkeller since
$${\sigma(2\cdot 5^2\cdot 7)-2(2\cdot 5^2\cdot 7) \over 2}=22=2\cdot 7+7+1$$ is a sum of factors of $2\cdot 5^2\cdot 7$ excluding $2\cdot 5^2\cdot 7$ and $5^2\cdot 7$; and $2\cdot 5\cdot 7^2$ is half-Zumkeller since 
$${\sigma(2\cdot 5\cdot 7^2)-2(2\cdot 5\cdot 7^2) \over 2}=23=2\cdot 7+7+2$$ is a sum of factors of $2\cdot 5\cdot 7^2$ excluding $2\cdot 5\cdot 7^2$ and $5\cdot 7^2$. \qed

The next four results are  pointed out in \cite{walsh}.

\begin{fact}\label{oddh}\cite{walsh} If $n$ is odd and half-Zumkeller, then $n$ is a perfect square.
\end{fact}

{\bf Proof.} If $n$ is odd and half-Zumkeller, then $\sigma(n)-n$ must be even and $\sigma(n)$ must be odd. Let the prime factorization of $n$ be $\Pi_{i=1}^m p_i^{k_i}$. Then 
$\sigma(n)=\Pi_{i=1}^m (\sum_{j=1}^{k_i} p_i^j)$. If $\sigma(n)$ is odd, then all $k_i$ must be even. So $n$ is a perfect square.
\qed

\begin{fact}\label{mnh}\cite{walsh} If $m$ and $n$ are half-Zumkeller numbers with $(m, n)=1$, then $mn$ is half-Zumkeller.
\end{fact}

\begin{fact} \label{h2}\cite{walsh}
If $n$ is even and half-Zumkeller, then $n$ is Zumkeller.
\end{fact}

{\bf Proof.}  Let $D$ be the set of all positive factors of $n$.
If $n$ is even and half-Zumkeller, then there exists $A\subset D\setminus\{n, {n \over 2}\}$ such that
$${n \over 2}+\sum_{a\in A} a=\sum_{b \in D, b\notin \{n, {n \over 2}\}\cup A} b.$$
Adding ${n \over 2}$ to both sides, we have
$$n+\sum_{a\in A} a={n \over 2}+\sum_{b \in D, b\notin \{n, {n \over 2}\}\cup A} b.$$
This means that all positive factors of $n$ are partitioned into two equal-summed parts and $n$ is Zumkeller. \qed

\begin{remark}\label{h4}\cite{walsh}
Let $n$ be even. $n$ is half-Zumkeller if and only if $n$ admits a Zumkeller parition such that $n$ and ${n \over 2}$ are in distinct parts.
\end{remark}

{\bf Proof.} Let $n$ be even. Let $D$ be the set of all positive factors of $n$.  $n$ is half-Zumkeller if and only if there exists $A\subset D\setminus\{n, {n \over 2}\}$ such that 
$${n \over 2}+\sum_{a\in A} a=\sum_{b \in D, b\notin \{n, {n \over 2}\}\cup A} b.$$
That is., $$n+\sum_{a\in A} a={n \over 2}+\sum_{b \in D, b\notin \{n, {n \over 2}\}\cup A} b.$$
This is equivalent to saying that $n$ admits a Zumkeller partition
such that $n$ and ${n \over 2}$ are in distinct parts. \qed

The following conjecture is proposed in \cite{walsh}.

\bigskip

{\bf Conjecture 2} \cite{walsh}
If $n$ is even and Zumkeller, then $n$ is half-Zumkeller.

\bigskip

In two of the next three results we shall verify that the conjecture is true in some cases. 

\begin{prop}\label{h5}
Let $n$ be even and Zumkeller. If $\sigma(n)<3n$, then $n$ is half-Zumkeller.
\end{prop}

{\bf Proof.} If $n$ is Zumkeller, by Fact \ref{f1}, ${\sigma(n)-2n \over 2}$ is either zero or the sum of some  factors of $n$ excluding $n$. If $\sigma(n)<3n$, then ${\sigma(n)-2n \over 2}<{n \over 2}$. So 
${\sigma(n)-2n \over 2}$ excludes ${n \over 2}$ as well. Then ${\sigma(n)-2n \over 2}$ is zero or the sum of some  factors of $n$ excluding $n$ and ${n \over 2}$. By Fact \ref{h3}, $n$ is half-Zumkeller. \qed

\begin{prop}\label{h6}
Let $n$ be even. Then, $n$ is Zumkeller if and only if $n$ is either  half-Zumkeller or ${\sigma(n)-3n \over 2}$ is $0$ or is a sum of some  factors of $n$ excluding $n$ and ${n \over 2}$.
\end{prop}

{\bf Proof.} Let $n$ be even. If $n$ is Zumkeller but not  half-Zumkeller , then by Remark \ref{h4}, any Zumkeller partition of the positive factors of $n$ must have $n$ and ${n \over 2}$ in the same parts. In other words, there exists $a$ which is either $0$ or a sum of some  factors of $n$ excluding $n$ and ${n \over 2}$ such that 
$$2(n+{n \over 2}+a)=\sigma(n).$$
So, $a={\sigma(n)-3n \over 2}$. Therefore, ${\sigma(n)-3n \over 2}$ is zero or a sum of some  factors of $n$ excluding $n$ and ${n \over 2}$.

If $n$ is half-Zumkeller, then $n$ is Zumkeller by Fact \ref{h2}. If ${\sigma(n)-3n \over 2}$ is zero or a sum of some  factors of $n$ excluding $n$ and ${n \over 2}$, then 
$${\sigma(n)-2n \over 2}={\sigma(n)-3n \over 2}+{n \over 2}$$ is a sum of some  factors of $n$ excluding $n$. By Fact \ref{f1}, $n$ is Zumkeller. \qed

\begin{prop}\label{addh5}
If 2 divides $n$, 3 divides $n$, $n$ is  Zumkeller, and $\sigma(n) < {10n \over 3}$, then $n$ is half-Zumkeller.
\end{prop}

{\bf Proof.}  If $n$ is not half-Zumkeller, by Proposition \ref{h6}, $(\sigma(n)-3n)/2$ is $0$ or a sum of some factors of $n$ excluding $n$ and $n/2$. Then, $$(\sigma(n)-2n)/2 = (\sigma(n)-3n)/2 + n/3 + n/6.$$ Since $\sigma(n)/n < 10/3$ we have that 
$(\sigma(n)-3n)/2 < n/6$. Hence $(\sigma(n)-2n)/2$  is a sum  of some factors of $n$ excluding $n$ and $n/2$. By Fact \ref{h3}, $n$ is half Zumkeller. This is a contradiction. \qed

The next proposition identifies some half-Zumkeller numbers.

\begin{prop}\label{twice}
If $n$ is Zumkeller, then $2n$ is half-Zumkeller.
\end{prop}

{\bf Proof.} Let  $n=2^kL$  with  $k$ a nonnegative integer and $L$ an odd number, be a Zumkeller number. Then all positive factors of $n$ can be partitioned into two disjoint equal-summed parts $D_1$ and $D_2$. Observe that every factor of $2n$ which is not a factor of $n$ can be written as $2^{k+1}l$ where $l$ is a factor of $L$. Observe that $2^kl$ is either in  $D_1$ or $D_2$. Without loss of generality, assume that $2^kl$ is in $D_1$. In this case, we move $2^kl$ to $D_2$ and add $2^{k+1}l$ to $D_1$. Perform this procedure to all factors of $2n$ which are not factors of $n$ except $2n$ itself. This procedure will yield an equal-summed partition of all factors of $2n$ except $2n$ itself. This shows that $2n$ is half-Zumkeller.
 \qed

The following remark is an immediate consequence of the above proposition.
It identifies the prime factorization of an even Zumkeller number
that is not half-Zumkeller.
\begin{remark}\label{remark54}
Let $n$ be even and the prime factorization of $n$ be $2^kp_1^{k_1}\cdots p_m^{k_m}$. If $n$ is Zumkeller but not half-Zumkeller, then  $2^ip_1^{k_1}\cdots p_m^{k_m}$ is not Zumkeller for any $i\le k-1$, and
$2^ip_1^{k_1}\cdots p_m^{k_m}$ is half-Zumkeller for any $i\ge k+1$. If Conjecture 2 is not true, then a counterexample to the Conjecture must have a prime factorization of $2^kp_1^{k_1}\cdots p_m^{k_m}$ such that $k$ is the minimum integer $i$ such that $2^ip_1^{k_1}\cdots p_m^{k_m}$ is Zumkeller.
\end{remark}

In order to investigate whether there are even Zumkeller numbers that are not half-Zumkeller, we shall see if all the results that are true for
Zumkeller numbers are true for even half-Zumkeller numbers also.

Results analogous to Propositions
\ref{fact55} and \ref{prop6}, are true for half-Zumkeller numbers also.

\begin{prop}\label{hprop55} Let $n$ be an even half-Zumkeller number and
$p$ be a prime with $(n, p)=1$. Then $np^l$ is half-Zumkeller for any
positive integer $l$.
\end{prop}

{\bf Proof.} Since $n$ is an even half-Zumkeller number, the set of all positive factors of $n$,
denoted by $D_0$ can be partitioned into two disjoint parts $A_0$ and $B_0$
so that the sums of  the two parts  are  equal and $n$ and ${n \over 2}$ are in distinct parts (by Remark \ref{h4}).

Group the factors of $np^l$ into $l+1$ groups $D_0$, $D_1$, $\ldots$, $D_{l}$  according to how
many factors of $p$ they admit, i.e., $D_i$ consists of all factors of $np^l$ admitting $i$
factors of $p$. Then each $D_i$ can be partitioned into two disjoint parts so that the sums of the two
parts  are equal and $np^i$ and ${np^i \over 2}$ are in distinct parts  according to the Zumkeller partition of the set $D_0$.
Therefore all positive factors of $np^l$ can be partitioned into two disjoint parts so that the sum of  these two parts   equal and $np^l$ and ${np^l \over 2}$ are in distinct parts. By Remark \ref{h4},  $np^l$ is half-Zumkeller.  \qed

The following is a direct Corollary of Proposition \ref{hprop55}.

\begin{coro}If $n$ is an even half-Zumkeller number and $m$ is a positive integer with $(n, m)=1$, then $nm$ is half-Zumkeller.
\end{coro}

But this result is not true for odd half-Zumkeller numbers.

\begin{remark}If $n$ is an odd half-Zumkeller number and
$m$ is a positive integer with $(m, n)=1$, then $mn$ need not be a
half-Zumkeller number. For example, $3^2\times 5^2$ is half-Zumkeller
since all its proper positive factors can be partitioned into two equal
parts: $75+9+5=45+25+15+3+1$. But $3^2\times 5^2\times 7$ is not
half-Zumkeller by Fact \ref{oddh}.
\end{remark}

We shall prove Proposition\ref{prop6}  for even half-Zumkeller numbers.

\begin{prop} \label{h30} Let $n$ be an even half-Zumkeller number and the prime factorization of  $n$ be $p_1^{k_1}p_2^{k_2}\cdots p_m^{k_m}$. Then for positive integers $l_1, \ldots, l_m$,
$$p_1^{k_1+l_1(k_1+1)}p_2^{k_2+l_2(k_2+1)}\cdots p_m^{k_m+l_m(k_m+1)}$$ is half-Zumkeller.
\end{prop}

{\bf Proof.} It is sufficient to show that $p_1^{k_1+l_1(k_1+1)}p_2^{k_2}\cdots p_m^{k_m}$ is half-Zumkeller
if $n=p_1^{k_1}p_2^{k_2}\cdots p_m^{k_m}$ is an even half-Zumkeller number.
Assume that $n=p_1^{k_1}p_2^{k_2}\cdots p_m^{k_m}$ is even and half-Zumkeller, then
the set of all positive factors of $n$, denoted by $D_0$ can be partitioned into two disjoint
parts $A_0$ and $B_0$ so that the sums of  the two parts  are  equal and $n$
and ${n \over 2}$ are in distinct parts (by Remark \ref{h4}). Note that the
factors of $p_1^{k_1+l_1(k_1+1)}p_2^{k_2}\cdots p_m^{k_m}$ can be partitioned
into $l_1+1$ disjoint groups $D_i$, $0\le i\le l_1$, where  elements in
$D_i$ are obtained by multiplying $p_1^{i(k_1+1)}$ with elements in $D_0$.
Using the partition $A_0, B_0$ of $D_0$ we can partition every $D_i$
into two disjoint parts $A_i$ and $B_i$ so that the sums of  the corresponding
parts are  equal  and $np_1^{i(k_1+1)}$ and ${np_1^{i(k_1+1)} \over 2}$
are in distinct parts. Therefore,  the set of all positive factors of
$p_1^{k_1+l_1(k_1+1)}p_2^{k_2}\cdots p_m^{k_m}$ can be partitioned into two
disjoint equal-summed parts   and $p_1^{k_1+l_1(k_1+1)}p_2^{k_2}\cdots p_m^{k_m}$ and ${p_1^{k_1+l_1(k_1+1)}p_2^{k_2}\cdots p_m^{k_m} \over 2}$ are in distinct parts. By Remark \ref{h4}, $p_1^{k_1+l_1(k_1+1)}p_2^{k_2}\cdots p_m^{k_m}$ is half-Zumkeller. \qed

We shall now proceed to prove Proposition \ref{proppraczu} and Theorem \ref{theopz} for half-Zumkeller numbers.

We need a definition. This is similar to the concept of a practical number, but is related to half-Zumkeller numbers as practical numbers are related to Zumkeller numbers.

\begin{defi} We say that a positive integer $n$ is a {\em quasi-practical} number
if every
positive integer $  \le  \sigma(n)-n$ can be written as a sum of
distinct positive factors of $n$ excluding $n$.
\end{defi}

Clearly every prime number is a quasi-practical number. Some simple quasi-practical numbers are given by the following propositon.

\begin{prop}\label{2kstrong}For every nonnegative integer $k$,  $2^k$ is quasi-practical.
\end{prop}

{\bf Proof.} Since every integer $m<\sigma(2^k)-2^k=2^{k}-1$ can be written as a sum of factors of $2^{k-1}$ (take the binomial expansion of $m$), $2^k$ is quasi-practical. \qed

By Fact \ref{h1}, if $n$ is quasi-practical and $\sigma(n)-n$ is
even, then $n$ is half-Zumkeller.

Our next aim is to show that every practical number is quasi-practical. In order
to show this we need a technical lemma.

\begin{lemma}\label{Lemma1} Let $p$ be a prime. Let $l, A_0, A_1,    ...,   A_l$ be nonnegative integers.
Every nonnegative $M \leq  A_0 + A_1p +    ...  + A_lp^l$ can be written
as  $C_0 + C_1p +    ...  + C_lp^l$  for some
$C_0 \leq A_0, C_1 \leq A_1,    ...,   C_l \leq A_l$ if and only if
the following $l$ conditions
$$A_0 + 1 \geq  p,$$ 
$$A_0 + A_1p+ 1 \geq  p^2,$$ 
$$......$$
$$A_0 + A_1p +    ...  + A_{l-1}p^{l-1}+1 \geq p^l$$ are satisfied.
\end{lemma}

{\bf Proof.}Necessity: If every nonnegative $M \leq  A_0 + A_1p +    ...  + A_{l}p^{l}$ can
be written as  $C_0 + C_1p +    ...  + C_{l}p^{l}$  for
some $C_0 \leq A_0, C_1 \leq A_1,    ...   C_{l} \leq A_{l}$,
then  for any $0 \le  i \le l-1, A_0 + A_1p +    ...  + A_{i}p^{i}+1\le A_0 + A_1p +    ...  + A_{l}p^{l}$ and this
can
be writtien as $C_0 + C_1p +    ...  + C_lp^l$  for some $C_0 \leq A_0, C_1 \leq A_1,    ...   C_l \leq A_l$.
Note that for some $ i+1 \le j \le l$, $C_{j} \ne 0$ since, otherwise $C_0 + C_1p +    ...  + C_ip^i <A_0 + A_1p +    ...  + A_{i}p^{i}+1$.
Therefore, $A_0 + A_1p +    ...  + A_{i}p^{i}+1\ge p^{j+1} \ge p^{i+1}$.

We thus have the necessity.

Sufficiency: We shall prove this by induction. For $ l =0$ there is nothing to be shown.

Let $l \ge 1$. Take a nonnegative integer $M \leq  A_0 + A_1p +    ...  +A_{l}p^{l}$.

Find the largest $C_l \le A_l$ such that $C_lp^l \le M $. If $C_l = A_l, $ then,
$M-A_lp^l \le A_0 + A+1p + A_2p^2 + \cdots +A_{l-1}p^{l-1}$ and the induction will take care of the rest.
If $C_l < A_l$, then, $C_lp^l \le M < (C_{l} + 1)p^l$. This implies that
$M-C_lp^l < p^l$. But, from the hypothesis, $p^l -1 \le    A_0 + A_1p +    ...  + A_{l-1}p^{l-1} $. Again the induction
will take care of the rest.

\qed

We are now ready to prove a result  about quasi-practical numbers that will help us
show that every practical number is quasi-practical.

\begin{prop}\label{strong} Let $n$ be a practical and quasi-practical number.
Let $p$ be a prime with $(n,p)=1$ and $l$ be a positive integer.
Then $np^l$ is quasi-practical if and only if  $p\le \sigma(n)+1$.
\end{prop}

{\bf Proof.}  Necessity: If $np^l$ is quasi-practical, then the positive
integer $p-1<\sigma(np^l)-np^l$  is a sum of factors of $np^l$ excluding $np^l$.
Since $p-1<p$, $p-1$ must be a sum of factors of $n$, then $p-1\le \sigma(n)$.
So $p\le \sigma(n)+1$.

Sufficiency:  Let  $M$ be a postive integer less than $\sigma(np^l)-np^l$. Then
\begin{eqnarray}\label{eqm}
M<\sigma(np^l)-np^l&=&\sigma(n)(\sum_{i=0}^l p^i)-np^l \nonumber \\
&=&(\sum_{i=0}^{l-1} \sigma(n)p^i) +  (\sigma(n)-n)p^l .
\end{eqnarray}
Since $\sigma(n)+1 \ge p$, we get that

$$\sigma(n) + \sigma(n)p+ 1 \geq  p^2,$$ 
$$......$$
$$\sigma(n) + \sigma(n) p +    ...  + \sigma(n) p^{l-1}+1 \geq p^l$$.

Since $n$ is practical every positive integer $\le \sigma(n)$ is a sum of factors of $n$.
Since $n$ is also  quasi-practical every positive integer $\le \sigma(n) -n $ is a sum of
factors of $n$ excluding $n$.
Applying Lemma \ref{Lemma1} to (\ref{eqm}),
$M$ can be written as $C_0 + C_1p +    ...  + C_lp^l$  for some
$C_0 \leq \sigma(n), C_1 \leq \sigma(n),    ..., C_{l-1} \leq \sigma(n),   C_l \leq \sigma(n)-n$
if $p\le \sigma(n)+1$. Since $n$ is a practical number and a quasi-practical number,
then $C_l$ is a sum of factors of $n$ excluding $n$ itself and
each $C_i$, $1\le i\le l-1$  is a sum of factors of $n$.
Therefore, $M$ is a sum of factors of $np^l$ excluding $np^l$ itself.
So $np^l$ is quasi-practical if  $p\le \sigma(n)+1$. \qed

Now, we are ready to derive some consequences of the above result.

\begin{prop}\label{pracstrong} Every practical number is quasi-practical.
\end{prop}

 {\bf Proof:} By Fact \ref{pract1} any practical number $n$ looks like 
$n = p_1^{k_1}p_2^{k_2}\cdots p_m^{k_m}$ and $p_1<p_2<\ldots<p_m$ where
$p_1 = 2$ and
for $1\le i\le m-1$, $p_{i+1}\le \sigma(p_1^{k_1}\cdots p_i^{k_i})+1$.  But $2^{k_1}$ is both practical
 and quasi-practical. Repeated use of Proposition \ref{strong} shows that n is also quasi-practical.

\qed

\begin{remark}\label{remarkpracqprac} A number is quasi-practical if and only if either
it is a practical number or a prime number. In fact even quasi-practical numbers are
same as practical numbers. The difficult part of this is Proposition \ref{pracstrong}.
We omit the rest of the details.

\end{remark}

We shall now find necessary and sufficient conditions for a practical number to be a half-Zumkeller
number.

\begin{prop}\label{propprachzu}
A practical number $n$ is a half-Zumkeller number if and only if $\sigma(n)$ is even.

\end{prop}
{\bf Proof.} If $n$ is practical and half-Zumkeller, it is an even half-Zumkeller number.
So, $n$ is Zumkeller by Fact \ref{h2}
and $\sigma(n)$ is even by Fact \ref{f3}.

If $\sigma(n)$ is even, then ${\sigma(n)-n \over 2}$ is an integer smaller
than $\sigma(n)-n$. Since $n$ is a practical number, it is quasi-practical
by Proposition \ref{pracstrong}.  Therefore
${\sigma(n)-n \over 2}$ is a sum of factors of $n$ excluding $n$. By
Fact \ref{h1}, $n$ is half-Zumkeller.

\begin{theo}\label{theophalf}
Let $n$ be a practical number, $l$ be a positive integer and $p$ be a
prime with $(n, p)=1$.

(i).  If $\sigma(n)$ is even, then $np^l$ is half-Zumkeller. 

(ii). If  $\sigma(n)$ is odd,   then $np^l$ is half-Zumkeller if and only if $p\le\sigma(n)$ and $l$ is odd.
\end{theo}

{\bf Proof.} (i). Let $n$ be a practical number.  Since $\sigma(n)$ is
even,  $n$ is half-Zumkeller by  Proposition \ref{propprachzu}.
By Proposition \ref{hprop55}, $np^l$ is half-Zumkeller.

(ii). $\sigma(n)$ is odd. If $np^l$ is half-Zumkeller, then $np^l$ is
Zumkeller by Fact \ref{h2}. By Theorem \ref{theopz}, $p\le\sigma(n)$
and $l$ is odd. If $p\le\sigma(n)$ and $l$ is odd, then $np^l$ is
practical by Fact \ref{pract1}. Since $l$ is odd,
$\sigma(np^l)=\sigma(n)\sum_{i=0}^l p^i$ is even.
By (i), $np^l$ is half-Zumkeller.
\qed

A result similar to Proposition \ref{multiplyz} can also be shown for
 half-Zumkeller numbers.

\begin{prop}\label{multiplyhz}
Let $n$ be  an even  integer and $p$ be a prime with $(n, p)=1$.  Then the following conditions are equivalent:

(i) $np$ is half-Zumkeller.

(ii) The set  of all positive factors of $n$ can be partitioned into two disjoint parts  $D_1$ and $D_2$ such that $n$ is in $D_1$, ${n \over 2}$ is in $D_2$, and  $p(\sum_{d\in D_2} d - \sum_{d\in D_1} d)$ can be written as a sum of some factors of $n$   minus  a sum of the rest of factors of $n$.

(iii) The set  of all positive factors of $n$ can be partitioned into two disjoint parts  $D_1$ and $D_2$ such that $n$ is in $D_1$, ${n \over 2}$ is in $D_2$, and ${(p+1)(\sum_{d\in D_2} d - \sum_{d\in D_1} d) \over 2}$ can be written as a sum of some elements in $D_2$ possibly minus  a sum of some elements in $D_1$.

(iv) The set  of all positive factors of $n$ can be partitioned into four disjoint parts  $A_1$, $A_2$,  $A_3$ and $A_4$ such that $n$ is in either $A_3$ or $A_4$, ${n \over 2}$ is in either $A_1$ or $A_2$, and $(p+1)\sum_{d\in A_1} d +(p-1)\sum_{d\in A_2} d=(p+1)\sum_{d\in A_3} d +(p-1)\sum_{d\in A_4} d$.
\end{prop}

{\bf Proof.}  By Remark \ref{h4}, $np$ is half-Zumkeller if and only
if the set $D_0$ of all positive factors of $n$ can be partitioned into
$D_0=D_1 \cup D_2$ and $D_0=D_3 \cup D_4$ such that
$$p(\sum_{d\in D_2} d - \sum_{d\in D_1} d)+(\sum_{d\in D_3} d - \sum_{d\in D_4} d)=0,$$ and 
$n$ is in $D_2$ and ${n \over 2}$ is in $D_1$. The rest of the proof
follows along the lines of the proof of Proposition \ref{multiplyz}. \qed

Analogous to Fact \ref{257} we have,

\begin{fact}\label{257h}
 $2\times 5^2\times 7^2\times p$ is half-Zumkeller for $p=11, 13$. 
\end{fact}

{\bf Proof.} In the proof that $2\times 5^2\times 7^2\times p$ is Zumkeller for $p=11, 13$. 
the partition of the set of all positive factors of
$n=2\times 5^2\times 7^2 = 2450 $ into $D_1$ and $D_2$ satisfies the
conditions that $n = 2450$ is in $D_1$ and ${n \over 2}=1225$ is in $D_2$ in
addition to condition (iii) in Proposition \ref{multiplyz}. So it satisfies
condition (iii) in Proposition \ref{multiplyhz} and
$2\times 5^2\times 7^2\times p$ is half-Zumkeller for
$p=11, 13$. 
\qed

\bigskip

We also have a result similar to Proposition \ref{prop7} for half-Zumkeller numbers.

\begin{prop}\label{hprop7}
If $1=a_1<a_2<\cdots<a_k=n$ are all factors of an  even number $n$ with $a_{i+1}\le 2a_i$ for all $i$ and $\sigma(n)$ is even, then $n$ is half-Zumkeller. 
\end{prop}

{\bf Proof.} Note that in the proof of Proposition \ref{prop7}, $a_k=n$ and $a_{k-1}={n \over 2}$ have different signs. So we get a Zumkeller partition of $n$ such that $a_k=n$ and $a_{k-1}={n \over 2}$ are in distinct parts. By remark \ref{h4}, $n$ is half-Zumkeller. \qed 

As in remark \ref{n!z}, using Proposition \ref{hprop7}, one can easily show that $n!$ is half-Zumkeller.

\bigskip

Our next result is similar to Remark \ref{remark54}. It gives more
necessary conditions for a number to be an even Zumkeller
number but not half-Zumkeller. 

\begin{prop}\label{znoth} Let $n$ be an even Zumkeller with the prime
factorization $2^kp_1^{k_1}p_2^{k_2}\cdots p_m^{k_m}$, where $2<p_1<\cdots<p_m$.
If $n$ is not half-Zumkeller, then there exists $i$, $1\le i\le m$ such
that $p_{i}>\sigma(2^k\cdots p_{i-1}^{k_{i-1}})+1$. If $j$ is  the
smallest such $i$, then  $k_1,\cdots, k_{j-1}$ must be even
and $j\le m-1$.
\end{prop}

{\bf Proof.} Since $n$ is Zumkeller, $\sigma(n)$ is even.
If $p_{i}\le \sigma(2^k\cdots p_{i-1}^{k_{i-1}})+1$ holds for every $i$,
$1\le i\le m$, then  $n$ is practical. By Theorem \ref{theophalf} (i),  $n$ is
half-Zumkeller. Therefore, there exists $i$, $1\le i\le m$ such that
$p_{i}> \sigma(2^k\cdots p_{i-1}^{k_{i-1}})+1$.  Let $j$ be smallest such $i$.
Then $2^kp_{1}^{k_{1}}\cdots p_{j-1}^{k_{j-1}}$ is a practical number.
If one of $k_1,\cdots, k_{j-1}$ is odd, then $\sigma(2^kp_{1}^{k_{1}}\cdots p_{j-1}^{k_{j-1}})$ is even.
By Theorem \ref{theophalf},  $2^kp_{1}^{k_{1}}\cdots p_{j-1}^{k_{j-1}}$ is
half-Zumkeller. By Proposition \ref{hprop55},  $2^kp_1^{k_1}p_2^{k_2}\cdots p_m^{k_m}$ is
half-Zumkeller. This contradicts the hypothesis. So $k_1,\cdots, k_{j-1}$ must be even.
If $j=m$, then $2^kp_{1}^{k_{1}}\cdots p_{m-1}^{k_{m-1}}$ is practical and
$\sigma(2^kp_{1}^{k_{1}}\cdots p_{m-1}^{k_{m-1}})$ is odd and $p_m>\sigma(2^kp_{1}^{k_{1}}\cdots p_{m-1}^{k_{m-1}})+1$,
then by Theorem \ref{theopz}(ii), $2^kp_1^{k_1}p_2^{k_2}\cdots p_m^{k_m}$ is
not Zumkeller and this contradicts the hypothesis that $n$ is Zumkeller. \qed

\bigskip

Next we show that Conjecture 2 is true for $n<7.233498900\times 10^9$. One can also show
that it is true for $n=7.233498900\times 10^9$. But we omit the details.

\begin{prop} Let $n$ be even and Zumkeller. If $n$ is not half-Zumkeller, then $n\ge 7.233498900\times 10^9$.
\end{prop}

{\bf Proof.} Let the prime factorization of $n$ be $2^kp_1^{k_1}p_2^{k_2}\cdots p_m^{k_m}$ with $2<p_1<\cdots<p_m$.
By Proposition \ref{znoth}, there exists $i$, $1\le i\le m$ such that
$p_{i}>\sigma(2^k\cdots p_{i-1}^{k_{i-1}})+1$. Let $j$ be the smallest such $i$.
Then  $k_1,\cdots, k_{j-1}$ must be even and $j\le m-1$. We will use this $j$ several times in the proof.

We will dicuss the cases $k=1, 2$ and $k\ge 3$ below.

Case 1. $k=1$. In this case, $n=2p_1^{k_1}p_2^{k_2}\cdots p_m^{k_m}$. 

Subcase 1(a): $m\ge 7$.

If $p_1\neq 3$, then $p_1\ge 5$.  

If $7\le m\le 10$, then  $p_1=5$ and $p_2=7$. Since 
if $p_1\ge 7$, then by Fact \ref{sigma},
$${\sigma(n) \over n}<(1+{1 \over 2})\prod_{i=1}^m {p_i \over p_i-1}<{3 \over 2}\cdot{7 \over 6}\cdot{11 \over 10}\cdot{13 \over 12}\cdot{17 \over 16}\cdot{19 \over 18}\cdot{23 \over 22}\cdot{29 \over 28}\cdot{31 \over 30}\cdot{37 \over 36}\cdot{41 \over 40}<3.$$
By Proposition \ref{h5}, $n$ is half-Zumkeller and this is a contradiction to our assumption. So $p_1=5$ when $m\le 10$. If $p_2>7$, then $p_2\ge 11$ and 
$${\sigma(n) \over n}<(1+{1 \over 2})\prod_{i=1}^m {p_i \over p_i-1}<{3 \over 2}\cdot{5 \over 4}\cdot{11 \over 10}\cdot{13 \over 12}\cdot{17 \over 16}\cdot{19 \over 18}\cdot{23 \over 22}\cdot{29 \over 28}\cdot{31 \over 30}\cdot{37 \over 36}\cdot{41 \over 40}<3.$$
By Proposition \ref{h5}, $n$ is half-Zumkeller and this is a contradiction to our assumption. So $p_2=7$ when $m\le 10$ and $p_1=5$. 

By Fact \ref{257h1},  $2\cdot 5\cdot 7$, $2\cdot 5^2\cdot 7$ and $2\cdot 5\cdot 7^2$ are  half-Zumkeller numbers.  By Proposition \ref{h30}, $2\cdot 5^{odd}\cdot 7^{odd}$, $2\cdot 5^{odd}\cdot 7^2$, and $2\cdot 5^{2}\cdot 7^{odd}$ are half-Zumkeller numbers. 
By Fact \ref{257h}, $2\cdot 5^2\cdot 7^2\cdot p$ is half-Zumkeller for $p\le 13$. Recall that $m\ge 7$, we have that 
$$n\ge \min\{2\cdot 5^2\cdot 7^2\cdot 17\cdot 19 \cdot 23\cdot 29\cdot 31,  2\cdot 5^4\cdot 7^2\cdot  11\cdot 13 \cdot 17\cdot 19\cdot 23\}>1.6\times 10^{10}.$$

If $m\ge 11$, then 
$n\ge 2\cdot 5\cdot 7\cdot 11 \cdot 13\cdot 17\cdot 19\cdot 23\cdot 29 \cdot 31\cdot 37\cdot 41>10^{14}$.

If $p_1=3$, then   $2\le j\le m-1$ (recall the meaning of $j$ defined in
the beginning of the proof) and recall that $k_1, \cdots, k_{j-1}$ must be even. So
\begin{eqnarray*}
n&>&2p_1^{k_1}\cdots p_{j-1}^{k_{j-1}}\sigma(2p_1^{k_1}\cdots p_{j-1}^{k_{j-1}})^{m-j+1} \\
&\ge& 2p_1^{2}\cdots p_{j-1}^{2}\sigma(2p_1^{2}\cdots p_{j-1}^{2})^{m-j+1} \\
&\ge& 2\cdot 3^{2(j-1)} \cdot [(1+2)(1+3+3^2)^{j-1}]^{m-j+1} \\
&\ge&2\cdot 3^{2(j-1)}(3\cdot 13^{j-1})^{m-j+1} \\
&=& 2\cdot 3^{m+j-1}\cdot 13^{(j-1)(m-j+1)} \\
&=& 2\cdot 3^{m+j-1+(j-1)(m-j+1)\log_3 13}.
\end{eqnarray*}
Since the exponent in the above expression is a quadratic function of $j$ with negative coefficient of $j$, it reaches the minimum at either $j=2$ or $j=m-1$ (recall that $2\le j\le m-1$). By comparing the function values at these two places, we get that it has the minimum value when $j=2$, so 
$$n>2\cdot 3^{m+1}\cdot 13^{m-1} \ge 2\cdot 3^{8}\cdot 13^{6}>6.3\times 10^{10}.$$

Subcase 1(b): $m\le 6$. 

If $p_1\ge 5$, then by Fact \ref{sigma},
$${\sigma(n) \over n}<(1+{1 \over 2})\prod_{i=1}^m (1+{1 \over p_i-1})<{3 \over 2}\cdot{5 \over 4}\cdot{7 \over 6}\cdot{11 \over 10}\cdot{13 \over 12}\cdot{17 \over 16}\cdot{19 \over 18}<3.$$ By Proposition \ref{h5}, $n$ is half-Zumkeller and this is a contradiction to our assumption. So 
$$p_1=3.$$
So $j\ge 2$. 
 Then by Proposition \ref{znoth}, $m\ge j+1\ge 3$.  By Proposition \ref{addh5},
$${\sigma(n) \over n}\ge {10 \over 3}.$$ 
If $p_2>\sigma(2\cdot p_1^{k_1})+1\ge \sigma(2\cdot 3^{2})+1=40$, then $p_2\ge 41$ and by Fact \ref{sigma}
$${\sigma(n) \over n}<(1+{1 \over 2})\prod_{i=1}^m {p_i \over p_i-1}<{3 \over 2}\cdot({41 \over 40})^5<{10 \over 3}.$$ So
$$p_2\le \sigma(2\cdot p_1^{k_1})+1.$$ 
So $j\ge 3$ and $m\ge j+1\ge 4$.

If $p_3>\sigma(2\cdot 3^{k_1}5^{k_2})+1\ge \sigma(2\cdot 3^{2}5^{2})+1=1210$, then 
$${\sigma(n) \over n}<{3 \over 2}\cdot{3 \over 2}\cdot{5 \over 4}\cdot({1211 \over 1210})^4<{10 \over 3}.$$ So
$p_3\le \sigma(2\cdot 3^{2}5^{2})+1$. Therefore $j\ge 4$ and  $m\ge j+1\ge 5$ by Proposition \ref{znoth}.
If $p_4>\sigma(2\cdot3^{k_1}5^{k_2}p_3^{k_3})+1\ge \sigma(2\cdot 3^{2}5^{2}7^2)+1=68914$, then 
$${\sigma(n) \over n}<{3 \over 2}\cdot{3 \over 2}\cdot{5 \over 4}\cdot{7 \over 6}\cdot({68915 \over 68914})^3<{10 \over 3}.$$ So
$p_4\le \sigma(2\cdot 3^{k_1}5^{k_2}p_3^{k_3})+1$. Therefore, $j\ge 5$ and  by Proposition \ref{znoth}, $m\ge j+1\ge 6$. Since $m=6$, then $j=5$ and   
$$p_5>\sigma(2\cdot3^{k_1}5^{k_2}p_3^{k_3}p_4^{k_4})+1\ge \sigma(2\cdot 3^{2}5^{2}7^{2}11^{2})+1=9165430.$$
So
$$n>2\cdot 3^{2}5^{2}7^{2}11^{2}(9165431)^2>2\times 10^{22}.$$

Case 2. $k=2$. In this case, $n=2^2p_1^{k_1}p_2^{k_2}\cdots p_m^{k_m}$. 

Subcase 2(a): $m\ge 7$.

If $p_1\le 7\le \sigma(2^2)+1$, then $2\le j\le m-1$. Recall that all $k_1$, $\cdots$, $k_{j-1}$ are even. So
\begin{eqnarray*}
n&>&2^2p_1^{k_1}\cdots p_{j-1}^{k_{j-1}}\sigma(2^2p_1^{k_1}\cdots p_{i-1}^{k_{j-1}})^{m-j+1} \\
&\ge& 2^2p_1^{2}\cdots p_{j-1}^{2}\sigma(2^2p_1^{2}\cdots p_{j-1}^{2})^{m-j+1} \\
&\ge&2^2\cdot 3^{2(j-1)}(7\cdot 13^{j-1})^{m-j+1} \\
&=& 2^2\cdot 3^{2j-2}\cdot 7^{m-j+1}\cdot 13^{(j-1)(m-j+1)} \\
&=& 2^2\cdot 3^{2j-2+(m-j+1)\log_3 7+(j-1)(m-j+1)\log_3 13}.
\end{eqnarray*}
Since the exponent in the above expression is a quadratic function of $i$ with negative coefficient of $i$, it reaches the minimum at either $j=2$ or $j=m-1$. By comparing the function values at these two places, we get that it has the minimum value when $j=2$, so 
$$n>2^2\cdot 3^2\cdot 7^{m-1}\cdot 13^{m-1} \ge 36\cdot 7^{6}\cdot 13^{6}>2\times 10^{13}.$$

If $p_1\ge 11$, then $m\ge 14$. Since if $m\le 13$ then 
$${\sigma(n) \over n}<(1+{1 \over 2}+{1 \over 4})\prod_{i=1}^m {p_i \over p_i-1}<{7 \over 4}\cdot{11 \over 10}\cdot{13 \over 12}\cdot{17 \over 16}\cdot{19 \over 18}\cdot{23 \over 22}\cdot{29 \over 28}\cdot{31 \over 30}\cdot{37 \over 36}\cdot{41 \over 40}\cdot{43 \over 42}\cdot{47 \over 46}\cdot{53 \over 52}\cdot{59 \over 58}<3.$$
By Proposition \ref{h5}, $n$ half-Zumkeller and this is a contradiction to our assumption. So $m\ge 14$ if $p_1\ge 11$. In this case,
$$n\ge 2^2\cdot 11 \cdot 13\cdot 17\cdot 19\cdot 23\cdot 29 \cdot 31\cdot 37\cdot 41\cdot 43\cdot 47\cdot 53\cdot 59\cdot 61>2.2\times 10^{21}.$$

Subcase 2(b): $m\le 6$. 

If $m\le 3$, then $p_1=3$ since otherwise, ${\sigma(n) \over n}<(1+{1 \over 2}+{1 \over 2^2})\cdot{5 \over 4}\cdot{7 \over 6}\cdot{11 \over 10}={7 \over 4}\cdot {5 \over 4}\cdot{7 \over 6}\cdot{11 \over 10} <3$ and $n$ is half-Zumkeller by Proposition \ref{h5}, this contradicts our assumption. So $j\ge 2$ and $m\ge j+1\ge 3$ by Proposition \ref{znoth}. Since $p_1=3$, by Proposition \ref{addh5},
${\sigma(n) \over n}\ge {10 \over 3}.$ 

If $m=3$: If $p_2> \sigma(2^2p_1^{k_1})+1\ge \sigma(2^23^2)+1=92$, then ${\sigma(n) \over n}<{7 \over 4}\cdot{3 \over 2}\cdot{97\over 96}\cdot {101 \over 100}<{10 \over 3}$.  Therefore, $p_2\le \sigma(2^2p_1^{k_1})+1$.  So $j\ge 3$ and $m\ge j+1\ge 4$.

What left is to consider $m=4$ or $5$ or $6$. If $p_1>\sigma(2^2)+1=8$, then $p_1\ge 11$ and
$${\sigma(n) \over n}<{7 \over 4}\cdot{11 \over 10}\cdot{13 \over 12}\cdot{17 \over 16}\cdot{19 \over 18}\cdot{23 \over 22}\cdot{29 \over 28}< 3.$$ 
So $p_1\le \sigma(2^2)+1$. Therefore, $2\le j\le m-1$. 

If $j=2$, then $p_2>\sigma(2^2p_1^{k_1})+1\ge \sigma(2^23^2)+1=92$, so $p_2\ge 97$  and 
$${\sigma(n) \over n}<{7 \over 4}\cdot{3 \over 2}\cdot{97 \over 96}\cdot{101 \over 100}\cdot{103 \over 102}\cdot{107 \over 106}< 3.$$ 
Therefore, $j\ge 3$.

If $j=3$, then $p_3>\sigma(2^2p_1^{k_1}p_2^{k_2})+1\ge \sigma(2^23^25^2)+1=2822$ and $m\ge j+1\ge 4$. So $p_3\ge 2833$, $p_4\ge 2837$, and
$$n\ge 2^23^25^2\cdot 2833\cdot 2837=7.202859300\times 10^9.$$
It can be shown that $2^23^25^2\cdot 2833\cdot 2837$ is half-Zumkeller. We omit the details.

If $j=4$, then $p_4>\sigma(2^2p_1^{k_1}p_2^{k_2}p_3^{k_3})+1\ge\sigma(2^23^25^27^2)+1=160798$ and $m\ge j+1\ge 5$. So 
$$n>2^23^25^27^2160799^2> 1.1\times 10^{15}.$$

If $j=5$ (only for $m=6$), then $p_5>\sigma(2^2p_1^{k_1}p_2^{k_2}p_3^{k_3}p_4^{k_4})+1\ge \sigma(2^23^25^27^211^2)+1=21386002$. So 
$$n>2^23^25^27^211^221386003^2> 2.4\times10^{21}.$$

Case 3. $k\ge 3$. In this case, $n=2^3p_1^{k_1}p_2^{k_2}\cdots p_m^{k_m}$.

Subcase 3a: $m\ge 7$.

 If $p_1>\sigma(2^k)+1\ge\sigma(2^3)+1=16$, then $p_1\ge 17$ and
$$n\ge 2^3\cdot 17\cdot 19\cdot 23\cdot 29\cdot 31\cdot 37\cdot 41>8.1\times 10^{10}.$$
If $p_1\le \sigma(2^k)+1$, then $2\le j\le m-1$. Recall that $k_1$, $\ldots$, $k_{j-1}$  must be even. 
So 
\begin{eqnarray*}
n&>&2^kp_1^{k_1}\cdots p_{j-1}^{k_{j-1}}(\sigma(2^kp_1^{k_1}\cdots p_{j-1}^{k_{j-1}}))^{m-j+1} \\
&\ge&2^33^{2j-2}(15\cdot 13^{j-1})^{m-j+1} \\
&\ge& 2^33^{2}15^{m-1}13^{m-1}\ge  2^33^{2}15^{6}13^{6}>3.9\times 10^{15}.
\end{eqnarray*}

Subcase 3b: $m\le 6$. 

If $p_1>7$, then $p_1\ge 11$ and 
$${\sigma(n) \over n}<2\cdot{11 \over 10}\cdot{13 \over 12}\cdot{17 \over 16}\cdot{19 \over 18}\cdot{23 \over 22}\cdot{29 \over 28}< 3.$$ 
So
$p_1\le7$. So $j\ge 2$ and  $m\ge j+1\ge 3$  by Proposition \ref{znoth}. 

We will discuss the cases $p_1=3$ or $5\le p_1\le 7$ separately as below. 

Subcase 3b1: $p_1=3$.

  If $p_2>\sigma(2^33^{k_1})+1\ge \sigma(2^33^2)+1=196$, then $p_2\ge 197$ and
$${\sigma(n) \over n}<2\cdot{3 \over 2}\cdot({197 \over 196})^5<{10 \over 3}.$$
By Proposition \ref{addh5}, $n$ is half-Zumkeller and this contradicts our assumption. So $p_2\le\sigma(2^33^{k_1})+1$. So $j\ge 3$ and  $m\ge j+1 \ge 4$  by Proposition \ref{znoth}.

So we can assume that $4\le m\le 6$. Recall that $3\le j\le m-1$ and $k_1, \cdots, k_{j-1}$ must be even.
So 
\begin{eqnarray*}
n&>&2^kp_1^{k_1}\cdots p_{j-1}^{k_{j-1}}(\sigma(2^kp_1^{k_1}\cdots p_{j-1}^{k_{j-1}}))^{m-j+1} \\
&\ge&2^33^2(5^{2})^{j-2}(15\cdot 13\cdot 31^{j-2})^{m-j+1} \\
&\ge& 2^33^{2}5^{2(m-3)}15^{2}13^{2}31^{2(m-3)} (j=m-1)\\
&\ge&  2^33^{2}5^{2}15^{2}13^{2}31^{2}>6.5\times 10^{10}.
\end{eqnarray*}

Subcase 3b2: $5\le p_1\le 7$.

 If $p_2>\sigma(2^3p_1^{k_1})+1\ge \sigma(2^35^2)+1=466$, then $p_2\ge 467$ and 
$${\sigma(n) \over n}<2\cdot{5 \over 4}\cdot({467 \over 466})^5< 3.$$
Since $n$ is not half-Zumkeller, then this is impossible by Proposition \ref{h5}. So $p_2\le\sigma(2^3p_1^{k_1})+1$. So $j\ge 3$  and $m\ge j+1\ge 4$ by Proposition \ref{znoth}. 

So we can assume that $4\le m\le 6$. Recall that $3\le j\le m-1$ and $k_1$, $\ldots$, $k_{i-1}$  must be even. 
So 
\begin{eqnarray*}
n&>&2^kp_1^{k_1}\cdots p_{j-1}^{k_{j-1}}(\sigma(2^kp_1^{k_1}\cdots p_{j-1}^{k_{j-1}}))^{m-j+1} \\
&\ge&2^35^{2j-2}(15\cdot 31^{j-1})^{m-j+1} \\
&\ge& 2^35^{4}15^{m-2}31^{2(m-2)} \ (j=3) \\
&\ge&  2^35^{4}15^{2}31^{4}>1.03\times 10^{12}. \qed
\end{eqnarray*}
 \section{Problems}

The following problems need further study.

1.\cite{walsh} Is every even Zumkeller number half-Zumkeller?

2. What are the odd Zumkeller numbers?

3. What are the odd half-Zumkeller numbers?

4. Does the set of Zumkeller numbers have density? Note that this set is a subset of abundant numbers and the set of abundant numbers has density between 0.2474 and 0.2480 \cite{marc}.


\begin{thebibliography}{JluR00}



\bibitem{walsh} S. Clark, J. Dalzell, J. Holliday, D. Leach, M. Liatti and M. Walsh, Zumkeller numbers, presented in the  Mathematical Abundance Conference at Illinois State Uiversity on April 18th, 2008.

\bibitem{marc} Marc Deleglise, Bounds for the density of abundant integers, Experimental Mathematics 7 (1998): 137-143.

\bibitem{rao} A. Galletti and K.P.S. Bhaskara Rao, A new proof of the Egyptian fraction theorem and great numbers, manuscript.

\bibitem{Srinivasan} A. K. Srinivasan, Practical numbers, Current Science 17 (1948): 179-180, MR0027799.


\bibitem{stewart} B. M. Stewart, Sums of distinct divisors, American Journal of Mathematics 76 (1954): 779-785.
   




\end{thebibliography}
\end{document}